\documentclass{amsart}
\usepackage{amsmath}
\usepackage{amsthm}
\usepackage{amscd}
\usepackage{amsfonts}
\usepackage{amssymb}
\usepackage[dvips,final]{graphics}
\newcommand{\IncludeXInvol}[1]{\includegraphics{XInvol#1.eps}}
\newcommand{\IncludeLoopPic}{\includegraphics{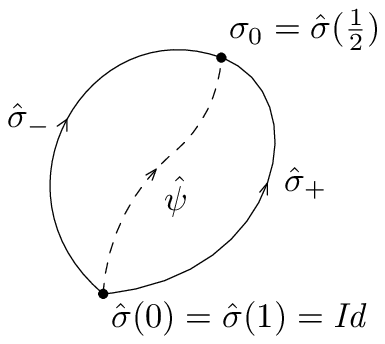}}

\newcommand{\iq}{\mathbf i}
\newcommand{\jq}{\mathbf j}

\newcommand{\hlinda}{{\mathbb H}}
\newcommand{\clinda}{{\mathbf{Ca}}}
\newcommand{\alinda}{{\mathcal{A}}}
\newcommand{\slinda}{{\mathcal{S}}}
\newcommand{\bl}{\mathrm{b}}

\newcommand{\rpn}{{\mathbb R} P^n}
\newcommand{\Z}{\mathbb{Z}}
\newcommand{\R}{\mathbb{R}}

\DeclareMathOperator{\Diff}{Diff}
\DeclareMathOperator{\rel}{rel}
\newcommand{\circdot}{\mathrel{\textstyle\circ\mkern-7.12mu\cdot}}

\newtheorem{theorem}{Theorem}
\newtheorem{lemma}{Lemma}
\newtheorem{proposition}{Proposition}
\theoremstyle{definition}
\newtheorem*{question}{Question}

\begin{document}

\title[Exotic involutions of Euclidean spheres]{Wiedersehen metrics and
exotic involutions of Euclidean spheres}
\author[Abresch]{U. Abresch$\, ^1$}
\address{Ruhr-Universit\"at Bochum,
Fakult\"at f\"ur Mathematik,
Universit\"atsstr.\ 150,
D-44780 Bochum, Germany}
\email{abresch@math.ruhr-uni-bochum.de}
\thanks{$^1$ Partially supported by DFG, FAEP, UNICAMP}
 
\author[Dur\'an]{C. Dur\'an$\, ^2$}
\address{IMECC-UNICAMP, Pra\c{c}a Sergio Buarque de Holanda, 651, 
Cidade Universit\'aria - Bar\~ao Geraldo, 
Caixa Postal: 6065 
13083-859 Campinas, SP, Brasil }
\email{cduran@ime.unicamp.br}
\thanks{$^2$ Supported by FAPESP grant 03/016789}

\author[P\"uttmann]{T. P\"uttmann$\, ^3$}
\address{Ruhr-Universit\"at Bochum,
Fakult\"at f\"ur Mathematik,
Universit\"atsstr.\ 150,
D-44780 Bochum, Germany}
\email{puttmann@math.ruhr-uni-bochum.de}
\thanks{$^3$ Partially supported by DFG, FAEP, UNICAMP}

\author[Rigas]{A. Rigas}
\address{IMECC-UNICAMP, Pra\c{c}a Sergio Buarque de Holanda, 651, 
Cidade Universit\'aria - Bar\~ao Geraldo, 
Caixa Postal: 6065 
13083-859 Campinas, SP, Brasil }
\email{rigas@ime.unicamp.br}
 

\begin{abstract}

We provide explicit, simple, geometric formulas for free involutions $\rho$
of Euclidean spheres that are not conjugate to the antipodal involution.
Therefore the quotient $S^n/\rho$ is a manifold that is homotopically
equivalent  but not diffeomorphic to $\mathbb R P^n$. We use these formulas
for constructing  explicit non-trivial elements in $\pi_1\mathrm{Diff}(S^5)$
and $\pi_1\mathrm{Diff}(S^{13})$ and to provide explicit formulas 
for non-cancellation phenomena in group actions.

\end{abstract}

\maketitle
 
\section{Introduction}

A smooth free involution $\rho$ on a sphere $S^n$ is called {\em exotic}
if it is not conjugate by a diffeomorphism to the standard antipodal
involution $\alpha(x) = -x$. The quotient $S^n/\rho$ is then a manifold
that is homotopically equivalent but not diffeomorphic to the standard
real projective space $\rpn$. 

There are several methods of constructing exotic involutions.
The first examples of such involutions were constructed by Hirsch and
Milnor (\cite{hirsch-milnor}) in $S^5$ and $S^6$, as restrictions to
invariant (standard) spheres of certain free involutions on exotic spheres.
Then there are examples constructed via surgery,  e.g. \cite{akbulut-kirby,
cappell-shaneson, fintushel-stern}. The reader can see also the basic
reference \cite{lopez-de-medrano} for topological and differentiable
invariants of involutions, and the classification and discussion using
analytical methods in \cite{atiyah-bott,oledzki}. 

A different path in the construction of exotic involutions is given
by simple involutions that restrict to involutions of Brieskorn
spheres (\cite{atiyah-bott, brieskorn, hirzebruch-mayer, hsiang}).
These have been used, for example, in the work of Grove and Ziller
\cite{grove-ziller} to construct metrics of non-negative sectional
curvature on exotic real projective spaces of dimension 5, and by Boyer,
Galicki, and Nakamaye \cite{boyer-galicki-nakamaye} to construct
Sasakian metrics of positive Ricci curvature on exotic real projective
spaces of dimension $4m+1$.

In this paper we construct free exotic involutions of Euclidean spheres
$S^n$, for $n = 5,6,13,14$. The origin of the formulas is quite geometric:
Recall that a Riemannian metric is called {\em wiedersehen} with respect
to points $N$ and $S$ in a manifold $M$ if every geodesic emanating from
$N$ reaches $S$ at a fixed length $\ell$ and vice versa. The wiedersehen
property at a point implies that $M$ is homeomorphic to the sphere
(see the book \cite{besse} for a complete discussion).
Lifting these geodesics to total spaces of bundles over spheres, one gets
many results and explicit formulas in differential and algebraic topology
(\cite{duran,duran-mendoza-rigas, classical,unitary}).

In dimensions 5 and 6, our involutions are essentially geometric formulas
for the Hirsch-Milnor involutions. These involutions are given by
restrictions of a natural involution of the Milnor exotic sphere
$\Sigma^7_{2,-1}$ to certain invariant submanifolds $\mathcal S_6$
and $\mathcal S_5$ that, by Morse theory, turn out to be spheres.
 
Using the geodesics of a wiedersehen metric of $\Sigma^7_{2,-1}$
constructed in \cite{duran}, we transfer the Hirsch-Milnor involutions
from the invariant 6-sphere $\mathcal S  $, now realized as the ``equator"
with respect to the metric, to the Euclidean 6-sphere contained in the
tangent space at a point of $\Sigma^7_{2,-1}$. The geometric origin of
these involutions is reflected in several features:

\begin{itemize}

\item
In contrast to the works above, our exotic involutions are described by
simple explicit formulas on the respective standard Euclidean spheres
(as opposed to Brieskorn spheres, or spheres inside of exotic spheres);
the simplicity of the formulas also translates to an elementary pictorial
description of the involutions (see Figures 1,2,3 in Section \ref{formulas}). 

\item
We can extend the constructions by substituting quaternions by Cayley
numbers everywhere, thus getting exotic involutions of spheres in
dimensions 13 and 14. These extensions are not immediate, since there
is no Cayley analog of the Gromoll-Meyer fibration
$Sp(2) \to \Sigma^7_{2,-1}$. This kind of phenomenon
-- the non-trivial extensions to the Cayley case --
already appears in \cite{duran-mendoza-rigas}.

\item 
The proof of the exoticity of the involutions is rather different from
the usual ones: for $S^6$ and $S^{14}$, we describe precisely the
$\mathbb Z_2$-action on $\Diff^+(S^n)$ of conjugation by the antipodal map.
This leads to several interesting open questions regarding the structure
of the relevant diffeomorphisms groups (see section \ref{sec-s6-s14}).
For $S^5$ and $S^{13}$, the involutions are shown to be exotic using
the fact that gluing diffeomorphisms $\sigma$ in $\pi_0\Diff^+(S^6)$
and $\pi_0\Diff^+(S^{14})$ have explicit lifts under the boundary map
$\pi_1 \Diff^+(S^n) \to \pi_0 \Diff^+(S^{n+1})$.
This is done in Section \ref{s5-s13}.

\item
The fact that these constructions admit lots of symmetries is exploited
in \cite{duran-puttmann}, where in particular we provide an explicit
cohomogeneity one diffeomorphism between a Brieskorn sphere and the
standard sphere $S^5$ which relates our constructions to the usual
Brieskorn ones. Again, extensions to the Cayley case are given there.
From these computations and known results
\cite{lopez-de-medrano,scharlemann-siebenmann} it follows that our
$S^5/\rho$ and $S^{13}/\rho$ are not homeomorphic to the standard
real projective spaces $\R P^5$ and $\R P^{13}$.

\item
Actually the discovery of the formulas for the involutions was somewhat
serendipitous, while we looked for geometric models of bundles over
exotic spheres. This is reflected by an application of these exotic
involutions, the construction of a very explicit example of non-cancellation
phenomenon in group actions: concretely, we will give non-conjugate actions
$r_1$, $r_2$ of $\mathbb Z_2 \times S^3$ on $X = S^6 \times S^3$
such that the restricted $\mathbb Z_2$- and $S^3$-actions are conjugate,
that is, ``neither factor can be cancelled"
(see Section \ref{noncancellation} for details).
Again, the formulas for these actions come from trivializations of
bundles using the geodesics of wiedersehen metrics.
\end{itemize}

\subsection*{Preliminaries}

We summarize here the fundamental topological facts that we need throughout the
paper (see e.g.~\cite{kervaire-milnor,kosinski}):
Let $\Diff^+(S^{n-1})$ and $\Diff^+(D^n)$ denote the group of orientation-preserving
diffeomorphisms of $S^{n-1}$ and $D^n$, respectively. Via the restriction
homomorphism, $\Diff^+(D^n)$ can be regarded as a normal subgroup of
$\Diff^+(S^{n-1})$. The quotient group
\begin{gather*}
  \Gamma_n = \Diff^+(S^{n-1})/\Diff^+(D^n) = \pi_0 \Diff^+(S^{n-1})
\end{gather*}
is abelian and consists of the equivalence classes of isotopic orientation
preserving diffeomorphisms of $S^{n-1}$.

The group $\Theta_n$ is the abelian group of h-cobordism classes of
homotopy $n$-spheres under the connected sum operation.
For $n\ge 5$, every homotopy $n$-sphere is homeomorphic to
$S^n$ and two homotopy $n$-spheres are h-cobordant if and only if
they are orientation preservingly diffeomorphic. Thus, $\Theta_n$
can be regarded as the group of all diffeomorphism classes of
differentiable structures on the topological $n$-sphere. For $n\ge 5$,
$\Theta_n$ is isomorphic to $\Gamma_n$. The isomorphism
$\Gamma_n \to \Theta_n$ is given by using $\sigma\in\Diff^+(S^n)$
to glue a twisted $n$-sphere from two disks.

Note that the group $\Theta_n$ (and thus the group $\Gamma_n$ for $n\ge 5$)
contains an important normal subgroup $bP^{n+1}$. This group consists of
all h-cobordism classes of homotopy $n$-spheres that bound
parallelizable manifolds.

\bigskip

\section{Explicit involutions of Euclidean spheres} \label{formulas}

Let $\hlinda$ and $\clinda$ denote the quaternions and the Cayley numbers,
respectively, and let $\Re$ and $\Im$ denote the real and imaginary parts.
Moreover, write 
\begin{gather*}
S^6 = \bigl\{ (p,w) \in \hlinda \times \hlinda 
  \,\, \big\vert \,\, \text{$\Re(p) = 0$, $|p|^2 + |w|^2 = 1$} \bigr\}\, ,
\end{gather*}
and similarly 
\begin{gather*}
S^{14} = \bigl\{ (p,w) \in \clinda \times \clinda
  \,\, \big\vert \,\, \text{$\Re(p) = 0$, $|p|^2 + |w|^2 = 1$} \bigr\}\, .
\end{gather*}

Consider the map $\bl: S^6 \to S^3$ (resp. $\bl:S^{14} \to S^7$) given by 
\begin{gather*}
\bl(p,w) = 
\begin{cases}  
 \frac {w}{|w|}e^{\pi p} \frac{\bar w}{|w|},&  w \neq 0\\
  -1, & w = 0,
\end{cases}
\end{gather*}
where $e^x$ denotes the exponential map of the group $S^3$ of unit
quaternions; thus $e^{\pi p} = \cos(\pi|p|) + \sin(\pi |p|) (p/|p|)$.
The map $\bl$ is a real analytic map whose homotopy class generates
$\pi_6 (S^3)$ (resp. $\pi_{14}(S^7)$, see \cite{duran-mendoza-rigas}).
We call these maps {\em Blakers-Massey elements}.
The map $\bl$ in the $6$-dimensional case is found using the wiedersehen
metric to explicitly represent the boundary map of the homotopy sequence
of the fibration $S^3 \cdots Sp(2) \to S^7$. This method has also been
used in \cite{classical, unitary} in order to produce explicit
representatives of several homotopy groups of the classical groups
along the borderline between the stable and unstable range.

Consider now $\sigma: S^6 \to S^6$ (resp. $S^{14} \to S^{14}$) given by
\begin{gather*}
\sigma(p,w) = (\bl(p,w) p \, \overline{\bl(p,w)}, \bl(p,w) w \, \overline{\bl(p,w)}) \, .
\end{gather*}
The map $\sigma$ is a real analytic, orientation-preserving diffeomorphism
that is {\em not} isotopic to the identity. Therefore the union of two
$7$-disks (resp $15$-disks) by $\sigma$ is an exotic sphere $\Sigma$.
This map is also found using the pointed wiedersehen metric
(in the $7$-dimensional case):
$\sigma = \exp_S^{-1}\circ \exp_N$,
where $N$ and $S$ denote the wiedersehen points of the metric
(see \cite{duran} for details). In the $7$-dimensional case,
$\Sigma$ generates the group $\Gamma_7 \cong \mathbb Z_{28}$,
and in the $15$-dimensional case, $\Sigma$ generates the first factor in
$\Gamma_{15} \cong  bP^{16} \times \mathbb Z_2
  \cong \mathbb Z_{8128} \times \mathbb Z_2$
(see \cite{duran,duran-mendoza-rigas}).

Let us consider now the map $\rho = \alpha\sigma$, where $\alpha$
is the antipodal map of the sphere. We have 

\begin{theorem} \label{rho-is-a-free-involution}
The map $\rho$ is a free involution of $S^6$ (resp. $S^{14}$).
\end{theorem}

\begin{proof}
Let us first write some remarkable properties of the map $\bl$ (see also
\cite{duran-mendoza-rigas}):

\begin{itemize}

\item
The map $\bl$ is equivariant under automorphisms of $\hlinda$
(resp. $\clinda$). In the quaternionic case this means that 
$\bl(qp\bar q, qw \bar q) = q \bl(p,w) \bar q$ for unit $q$.
In the case of Cayley numbers, the latter property holds provided
that $q$ lies in the subalgebra generated by $p$ and $w$. 

\item
$\bl(-p,-w) = \overline{\bl(p,w)}$.

\item
$\sigma^k (p,w) = \bigl(\bl(p,w)^k p \, \overline{\bl(p,w)^k},
\bl(p,w)^k w \, \overline{\bl(p,w)^k}\bigr)$.

\end{itemize}

With these properties the fact that $\rho\circ \rho$ is the identity
is easy to establish. 

In order to prove that $\rho$ is free, first note that $\rho(p,w) = (p,w)$
means that both $p$ and $w$ anticommute with $\bl(p,w)$. In particular
this implies that the real parts of $w$ and $\bl(p,w)$ are zero
(and the real part of $p$ is zero by definition).
Note that $e^{\pi p} = \cos(\pi |p|) + \sin(\pi |p|)\frac{p}{|p|}$;
thus ${\rm Re}(\bl(p,w)) = 0$ if and only if $|p|=1/2$ and thus
$|w| = \sqrt{3}/2$. Using the substitution
$p \mapsto qp\bar q$, $w \mapsto qw\bar q$,
we can assume that $p = \iq/2$ and
$w= \frac{{\sqrt{3}}}{2}(\cos(\theta)\iq + \sin(\theta)\jq)$
for some $\theta \in [0,\pi]$.
Then $\bl(p,w) = \cos(2\theta)\iq + \sin(2\theta)\jq$.
From the fact that $\bl(p,w)$ anticommutes with $p$ we get
$\cos(2\theta)=0$ and hence $\sin(2\theta)= \pm 1$.
Thus $\theta = \pi/4$ or $\theta = 3\pi/4$.
In none of these two case $w$ anticommutes with $\bl(p,w)$.
This argument generalizes to the Cayley case since we are dealing
with the algebra generated by $p$ and $w$ and therefore everything
happens in a copy of $\hlinda$ inside of $\clinda$.
\end{proof}

Note that the map $\rho$ restricts to the $5$-sphere (resp. $13$-sphere)
given by the condition $\Re(w) = 0$. In the next two sections we will
show that all of these involutions are exotic; we finish this section
by giving a pictorial sequence in Figures\,1~to~3 showing the involution
on $S^5$, done by translating the fact that the conjugation
$x \mapsto qxq^{-1}$ by a unit quaternion $q$ acting on a purely
imaginary $x$ is given by rotating $x$ along the axis $\Im(q)$
with angle $\theta$, where $\cos(\theta) = 2\Re(q)^2-1$.

\medskip

\begin{figure}
\begin{minipage}[b]{\linewidth}
\begin{center}
    \mbox{\scalebox{.4}{\IncludeXInvol{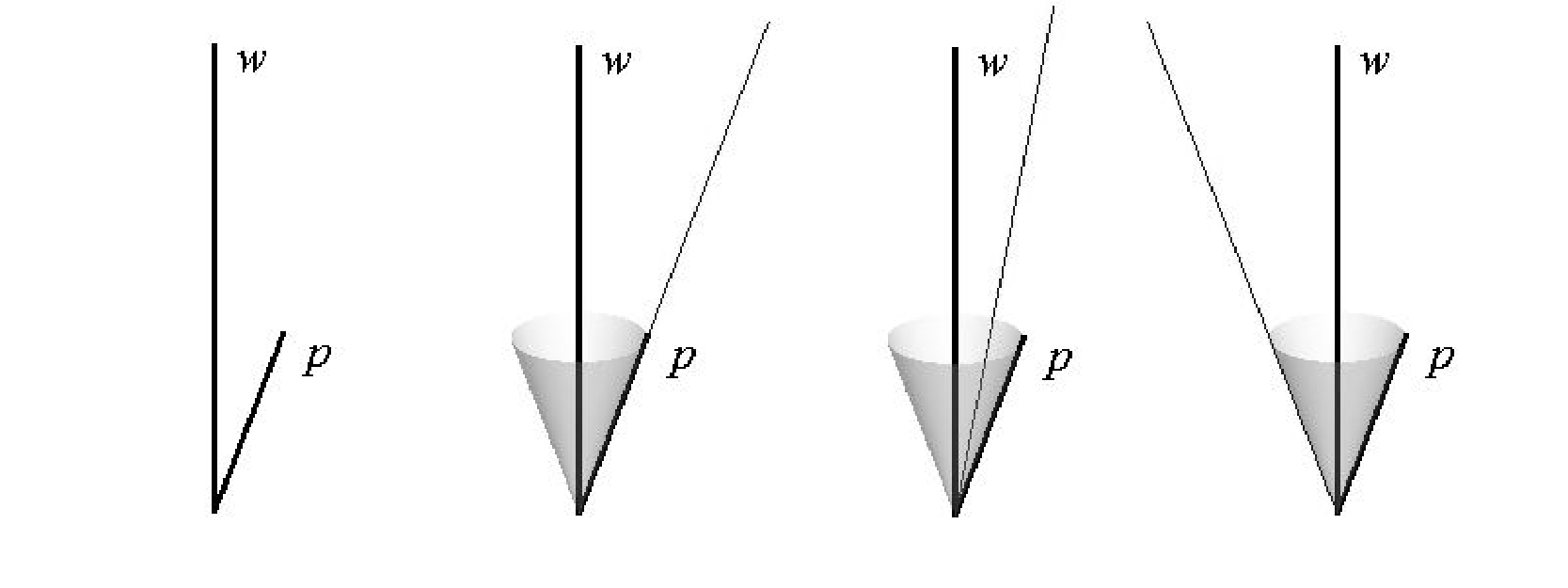}}}
\end{center}
\caption{Consider two vectors $p,w\in\R^3$ with $|p|^2 + |w|^2 = 1$.
The set of all such vectors forms the sphere $S^5$. Suppose that $w \neq 0$.
First rotate $p$ $180^{\circ}$ around $w$ and obtain an oriented axis.}
\end{minipage}
\begin{minipage}[b]{\linewidth}
\bigskip
\begin{center}
    \mbox{\scalebox{.4}{\IncludeXInvol{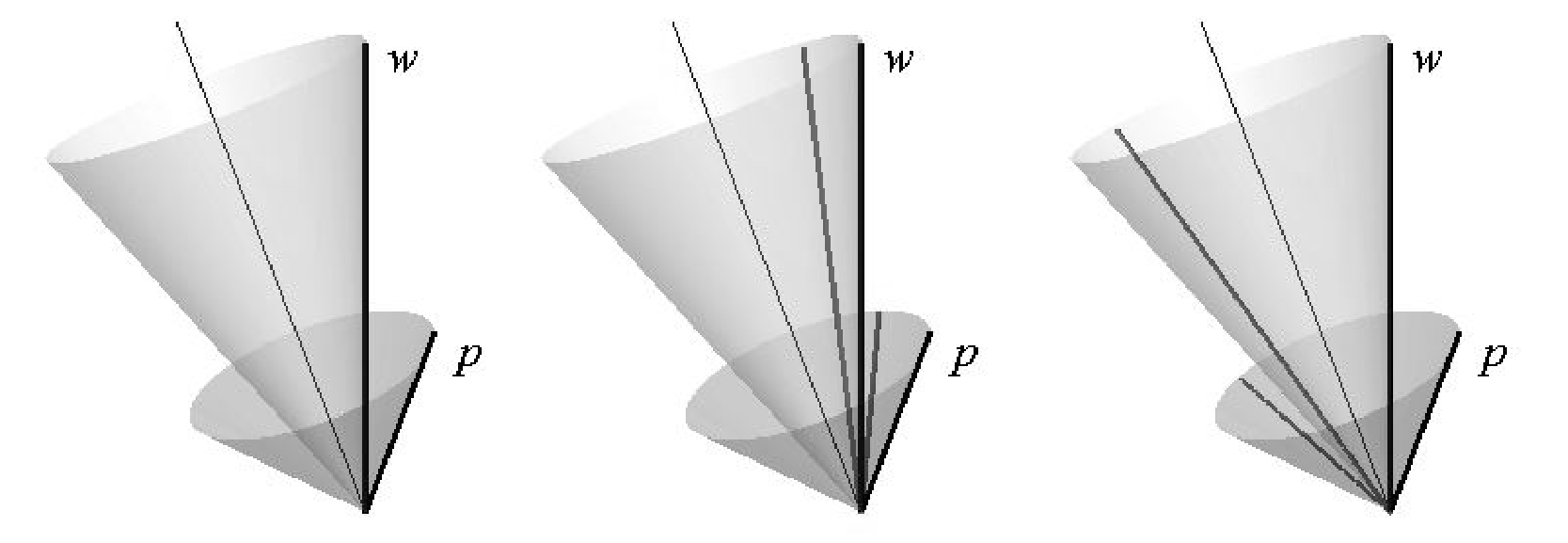}}}
\end{center}
\caption{Second, rotate $p$ and $w$ around this new axis by an angle of
$|p| \cdot 360^{\circ}$.}
\end{minipage}
\medskip
\begin{center}
    \mbox{\scalebox{.4}{\IncludeXInvol{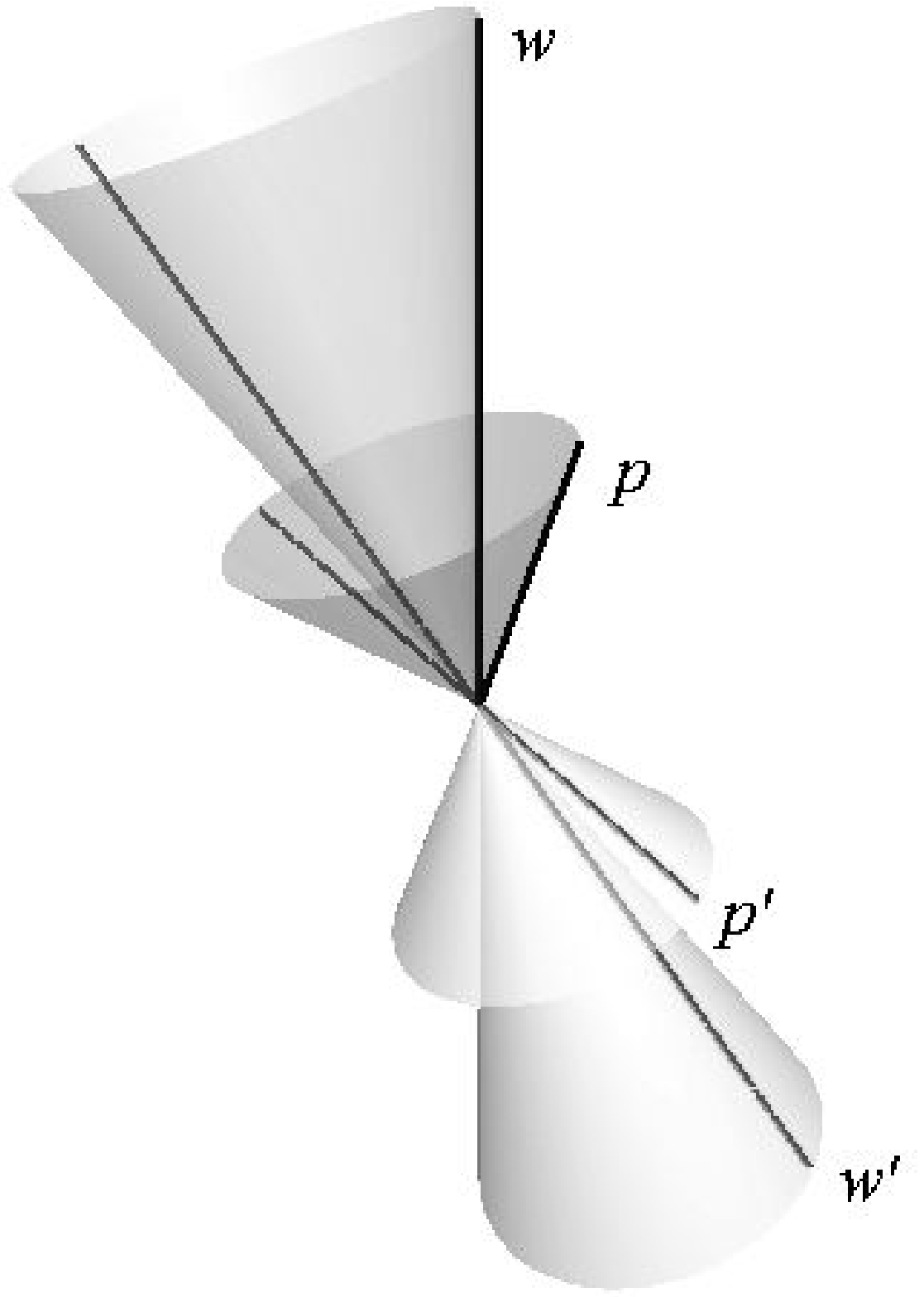}}}
\end{center}
\caption{Finally map the resulting vectors to their antipodes.
This process extends analytically to the case where $w=0$
(in this case $w' = 0$, $p' = -p$) and gives a visual description
of an exotic involution of $S^5$.}
\end{figure}

It can be read easily from Figures\,1~to~3 that the map constructed
is a fixed point free involution of $S^5$:
\begin{itemize}
\item
If $w \to 0$ then $|p| \to 1$. Hence in the second step we get close to a
rotation by $360^\circ$, i.e., to the identity of $S^5$. The map therefore
extends at least continuously to the case where $w=0$.
\item
If $w$ and $p$ are linearly dependent, the image of $w$ and $p$ is just
$-w$ and $-p$. If $w$ and $p$ are linearly independent the axis of
rotation constructed in the first step is still contained in the plane
spanned by $w$ and $p$.
Hence, they cannot be mapped to $-w$ and $-p$ in the second step.
This shows that the map is fixed point free.
\item
When the map is applied to $w'$ and $p'$ the axis constructed
in the first step is given by the central direction of the cones
opening downwards. In the second step one rotates around this axis
by an angle of $|p'|\cdot 360^{\circ}$ (observe the right hand rule).
Since $|p'| = |p|$ this gives $-w$ and $-p$.
Hence in the final step one is back at $w$ and $p$.
This shows that the map is an involution.
\end{itemize}

Of course, it would be nice to see from the pictorial description of
the involution why it is exotic
(without going up in dimension to $\Sigma^7_{2,-1}$).

\bigskip

\section{An involution of $\Diff^+(S^{n-1})$} \label{sec-s6-s14}

Now let us prove that $\rho$ is not equivalent (in the conjugation sense)
to the antipodal involution $\alpha$. In order to accomplish that,
consider the following $\mathbb Z_2$-action $\alinda$ on $\Diff^+(S^n)$, 
\begin{gather*}
\alinda (h) = \alpha \circ h \circ \alpha^{-1}
  = \alpha \circ h \circ \alpha \, ,
\end{gather*}
or, more simply, 
\begin{gather*}
\alinda (h)(x) = -h(-x) \, .
\end{gather*}

Note that, being conjugation, $\alinda$ is a group isomorphism;
$\alinda(hf) = \alinda(h)\alinda(f)$ in the group $\Diff^+(S^{n-1})$.
Note also that if two diffeomorphisms $h_0$ and $h_1$ are joined
by a curve of diffeomorphism $h_t$, then $\alinda(h_t)$ joins
$\alinda(h_0)$ to $\alinda(h_1)$ through diffeomorphisms.
Therefore $\alinda$ descends to an action on
$\pi_0 \Diff^+(S^{n-1}) = \Gamma_n$;
which by abuse of notation we also denote by $\alinda$.

In the particular case of $n = 7$ or $n=15$, we can explicitly compute
how this action behaves in $\Gamma_7 \cong \mathbb Z_{28}$ and in the
index $2$ subgroup $bP^{16}\cong\mathbb Z_{8128}$ of $\Gamma_{15}$
by acting on the known representative $\sigma$. In fact, Brumfiel
\cite{brumfiel} has shown that
$\Gamma_{15}\cong\Theta_{15}\cong
  bP^{16} \oplus \bigl(\Theta_{15}/bP^{16}\bigr) \cong Z_{8128}\oplus\Z_2$.
A computation shows the following:

\medskip

\noindent{\it Main commutation relation}:
$\alpha\sigma = \sigma^{-1}\alpha$, or equivalently, 
\begin{gather*}
\alinda(\sigma) = \sigma^{-1} \, .
\end{gather*}

\medskip

Thus the action on $\Gamma_7$ and on the subgroup $bP^{16}\subset \Gamma_{15}$
is given by 
\begin{gather*}
\alinda (n) = - n \, .
\end{gather*}
It would be very interesting to find how $\alinda$ acts on the
complement of $bP^{16} \subset \Gamma_{15}$. Clearly, it maps the elements
of order $2$ into themselves. So the question boils down to whether
$\alinda$ interchanges $(0,1)$ and $(4064,1)$ or not.
Alas, we do not have any explicit representative to decide this.

We are now ready to show that $\rho$ is not conjugate to the identity:

\begin{theorem} \label{rho-is-exotic}

The map $\rho$ is not conjugate to the antipodal map, i.e., 
there is no diffeomorphism $h:S^6 \to S^6$ (resp $h:S^{14} \to S^{14})$
that satisfies $h\rho h^{-1} = \alpha$.

\end{theorem}

\begin{proof}
 
Let us assume that such a diffeomorphism exists. Without loss of generality,
we can suppose that $h$ is orientation preserving, since if there exists an
orientation-reversing diffeomorphism $j$ such that $j\rho j^{-1} = \alpha$,
the diffeomorphism $h = \alpha j$ satisfies the same equation and is
orientation preserving. 
 
Now, we have,
\begin{gather*}
h \rho h^{-1}= \alpha \Leftrightarrow \sigma
  = \alpha h^{-1} \alpha h \Leftrightarrow \sigma
  = \alinda(h^{-1}) h \, .
\end{gather*}
 
Taking isotopy classes on both sides, we have 
\begin{gather*}
[\sigma] = [ \alinda(h^{-1}) h ]
  = [ \alinda(h^{-1})][ h ] = 2[h] + \tau \, ,
\end{gather*}
where the term $\tau$ only appears in the $14$-dimensional situation;
$\tau$ is zero except in the case that
$[h] = (n,1) \in  Z_{8128}\oplus\Z_2 \cong \Gamma_{15}$,
and $\alinda$ interchanges $(0,1)$ and $(4064,1)$. 
Then, $\tau = (4064,0)$. In any case, $\sigma$ is even inside of
$bP_8 \cong \mathbb Z_{28}$ (resp. $bP^{16} \cong \mathbb Z_{8128}$),
which contradicts the fact that $\sigma$ is a generator.
\end{proof}
 
Let us give some remarks on the proof. Notice that, with the same
computations as in Proposition \ref{rho-is-a-free-involution},
the maps $\rho_k = \alpha\sigma^k$ are also free involutions.
However, the proof of Theorem \ref{rho-is-exotic} only shows
that $\rho_k$ is exotic for odd $k$. Indeed this must be so,
since we have
\begin{gather*}
 \sigma^{-\ell} \rho_k \sigma^\ell =  \sigma^{-\ell} \alpha\sigma^k \sigma^\ell = \alpha\sigma^{k+2\ell} = \rho_{k+2\ell} \, ,
\end{gather*}
and therefore all the $\rho_k$ for $k$ even are conjugate among
themselves (in particular conjugate to the antipodal map $a$),
and all the odd $\rho_k$ are conjugate to $\rho_1 = a\sigma$. 
 
In fact, what we have done amounts to an analysis of the combinatorics
of the path-connected components of $\Diff(S^6)$ (and half of the
components of $\Diff(S^{14})$, but for simplicity we just discuss the
$S^6$-case). Indeed, $\Diff(S^6)$ has $56$ connected components, half of
which are orientation preserving and the other half orientation reversing.
The orientation preserving components are represented by the classes
$[Id], [\sigma], [\sigma^2], \dots, [\sigma^{27}]$
whereas the orientation reversing components are represented by the classes
$[\alpha], [\alpha\sigma], [\alpha\sigma^2], \dots , [\alpha\sigma^{27}]$.
We have that conjugation by the antipodal map acts, on each of these halves,
by $n \mapsto -n$ in $\mathbb Z_{28}$, that is, by fixing $[Id]$
(resp. [$\alpha$]) and the component of $\sigma^{14}$
(resp. [$\alpha\sigma^{14}]$), and permuting symmetrically the rest
$[\sigma^n] \mapsto [\sigma^{28-n}]$
(resp.  $[\alpha\sigma^n] \mapsto [\alpha\sigma^{28-n}])$.
As an additional by-product of this proof, we also get 
(compare \cite{mann-miller})

\begin{theorem}
 
Every orientation reversing diffeomorphism of $S^6$ is isotopic to a free involution.

\end{theorem}
 
Let us close this section with several remarks and questions from the global analysis
point of view:  The identity homeomorphism is a fixed point of the
$\mathbb Z_2$-action $\alinda$. Any orientation preserving fixed point
of this action must lie either on the component of the identity
or in the component of $\sigma^{14}$. Note that
$\alinda(\sigma^{14}) = \sigma^{-14}$, which is not the same map;
it just lies in the same path connected component. 

\begin{question}
Is there an odd map $f\in \Diff^+(S^6)$ in the isotopy class of $\sigma^{14}$?
That is, a fixed point of $\alinda$.  
\end{question}
 
\begin{question} 
Is there $f\in \Diff^+(S^6)$ in the isotopy class of $\sigma^{14}$ that satisfies
$f\circ f = Id?$. That is,  a fixed point of the inverse involution
$\mathcal B: \Diff^+(S^6) \to  \Diff^+(S^6) $ given by $\mathcal B(f) = f^{-1}$. 
Again, such a fixed point can only be isotopic either to the identity, or to $\sigma^{14}$. 
\end{question}
 
The relevance of such an involution is that it would be a diffeomorphism
that realizes the isotopy $\sigma^{28} \cong Id$ ``on the nose",
thus greatly helping in the understanding of exotic diffeomorphisms
and exotic spheres.

We can broaden these questions as follows:

\begin{question}
Which isotopy classes of orientation-preserving
diffeomorphisms can be realized by maps of finite order? 
\end{question}

For example, find a diffeomorphism $\eta:S^6 \to S^6$ representing
$4 \in \mathbb Z_{28} \cong \Gamma_7$ such that $\eta^7 = Id$. 

Explicit answers to these questions would provide exotic diffeomorphisms that  
improve upon the diffeomorphism $\sigma$, since they would express the group
structure of $\Gamma_7 \cong \mathbb Z_{28}$ in a direct way.

Also, note that the main commutation relation can be expressed as the statement
that the powers $\sigma^k$ of $\sigma$ are contained in the subset of
$\Diff^+(S^6)$ where the $\mathcal A$-orbit and the $\mathcal B$-orbit coincide. 
It would be interesting to study the structure of this subset; in particular, to find 
some other elements.

\bigskip

\section{Restriction to invariant spheres} \label{s5-s13}
 
As we have remarked, the exotic involution $\rho$ of $S^6$
(resp. $S^{14}$) has an invariant $5$-sphere (resp.\ invariant $13$-sphere)
given by the $\Re(w) = 0$. We have

\begin{theorem} \label{5-y-13}

The restriction of $\rho$ to $S^5$ (resp. $S^{13}$) is a free involution
that is not conjugate to the antipodal map.

\end{theorem}
 
The proof for $S^5$ could be copied from \cite{hirsch-milnor} in the
$5$-dimensional case. However, their proof breaks down in $S^{13}$ since
there are exotic 13-spheres. In the rest of this section we prove both cases
at once using the boundary map $\partial$ in the classical exact sequence
\begin{gather*}
  \pi_1 \Diff^+ (S^{n-1}) \stackrel{\partial}{\longrightarrow}
  \pi_0 \Diff^+ (S^n)
  \to \pi_0\Diff^+ (D^n) \to \pi_0\Diff^+ (S^{n-1}) \to \Gamma_n \to 0.
\end{gather*}
By a theorem of Cerf~\cite{cerf}, $\pi_0 \Diff^+ (D^n) = 0$ for $n\ge 6$.
In our context, Cerf's result implies that there are loops $\hat\sigma$
in $\Diff^+(S^5)$ and $\Diff^+(S^{13})$ which map under $\partial$
to diffeomorphisms of $S^6$ and $S^{14}$ that are isotopic to our
exotic diffeomorphisms $\sigma$. It is not a priori clear that
such loops $\hat \sigma$ can be given explicitly.

However, the concrete formula for $\sigma$ is of such a form,
as we shall see in Lemma\,\ref{explicit-loop} below.
In line with the spirit of this paper, we shall use these explicit loops
$\hat\sigma$ to prove Theorem\,\ref{5-y-13}; they give identities of maps
at several instances that induce very concrete homotopy identities.
Moreover, our proof becomes independent of Cerf's theorem this way.

Let us now recall a concrete definition of the boundary
homomorphism $\partial$. The elements of $\pi_1 \Diff^+(S^{n-1})$
can be represented by paths
$\hat \beta: [0,1] \to \Diff^+(S^{n-1})$
that map a neighborhood of $\{0,1\}$ to the identity map.
By standard approximation results, it can be assumed that such a path
$\hat \beta$ induces a smooth map $S^{n-1} \times [0,1] \to S^{n-1}$,
which is, by abuse of notation, again denoted by $\hat \beta$.
The image of $\hat\beta$ under the boundary map $\partial$
is now given by
\begin{align*}
  \partial(\hat \beta) : S^{n-1} \times [0,1] &\to S^{n-1} \times [0,1] \\
  (x,t) &\mapsto (\hat\beta(x,t),t) \, .
\end{align*}
Clearly, $\partial(\hat\beta)$ is an orientation preserving
diffeomorphism of the cylinder $S^{n-1} \times [0,1]$ which coincides
with the identity in a neighborhood of the boundary.
Thus $\partial(\hat \beta)$ induces an orientation preserving
diffeomorphism $\partial(\hat \beta)$ of the sphere $S^n$.

\begin{lemma}\label{explicit-loop}

The exotic diffeomorphisms $\sigma$ of $S^6$ and $S^{14}$
naturally define explicit loops $\hat \sigma$ in $\Diff^+(S^5)$
and $\Diff^+(S^{13})$ with $\partial(\hat\sigma) = \sigma$.

\end{lemma}

\begin{proof}

Recall that for $(p,w)$ satisfying $\Re(p) = 0$, one has 
\begin{gather*}
\sigma(p,w) = (\bl(p,w) p \, \overline{\bl(p,w)}, \bl(p,w) w \, \overline{\bl(p,w)}) \, .
\end{gather*}
Separating $w = w_0 + \omega$, where $w_0 = \Re(w)$, we have that 
\begin{gather*}
\sigma(p,w_0,\omega)
  = \bigl( \bl(p,w_0 + \omega) p \, \overline{\bl(p,w_0 + \omega)},
  \, w_0, \,
  \bl(p,w_0 + \omega) \omega \, \overline{\bl(p,w_0 + \omega)} \bigr) \, ,
\end{gather*}
since conjugation preserves the real part.
Thus we can think of $w_0\in [-1,1]$ as the parameter of a curve
of diffeomorphisms of the standard $S^5$ (resp.\ $S^{13}$).
Note that for $w_0 = \pm 1$, $\bl(p,w) = 1$, and thus at these levels
the diffeomorphisms are the identity.
\end{proof}

Let $T_n$ denote the characteristic subgroup in
$\pi_0 \Diff^+(S^n) \cong \Gamma_{n+1}$
generated by all elements of order $2$, and let
$N_{n-1} := \partial^{-1}(T_n)$
be the corresponding normal subgroup in $\pi_1 \Diff^+(S^{n-1})$.

\begin{lemma}\label{loop-two}

The loops $\hat\sigma$ introduced in Lemma\,\ref{explicit-loop}
generate the cyclic subgroups
\begin{gather*}
  \pi_1\Diff^+(S^5)/N_5 = \Z_{14}
  \quad \text{and}\quad
\pi_1 \Diff^+(S^{13})/N_{13} = \Z_{4064}.
\end{gather*}
\end{lemma}

\begin{proof}
Because of the definition of $N_{n-1}$, the boundary maps $\partial$
induce isomorphisms
\begin{gather*}
  \bar\partial: \pi_1\Diff^+(S^{n-1})/N_{n-1} \to \pi_0 \Diff^+(S^n)/T_n.
\end{gather*}
Now observe that by definition $T_6$ is the subgroup
$\Z_2 = \langle 14\rangle$ in $\pi_0 \Diff^+(S^{6})$
and that $T_{14}$ is the subgroup
$\Z_2\times\Z_2 = \langle (4064,0),(0,1)\rangle$ in
$\pi_0 \Diff^+(S^{14}) = bP^{16} \times \Z_2 = \Z_{8128} \times \Z_2$.
To conclude the argument, recall that the diffeomorphisms
$\sigma = \partial (\hat\sigma)$ generate $\Z_{28}$ and $bP^{16}$,
respectively.
\end{proof}

Before we employ this structural information in the proof of
Theorem\,\ref{5-y-13} we need to introduce some notation:
Let $\sigma_0$ be the restriction of $\sigma$ to the equator
of $S^6$ (resp.\ $S^{14}$) given by $w_0 = \Re(w) = 0$ and
let $\alpha_0$ denote the antipodal map of the equator.
In order to be consistent with the standard convention for
concatenation of loops we will now define all loops
and paths on the unit interval $[0,1]$. In particular,
we assume that $\hat \sigma$ is parametrized on $[0,1]$ such that $\hat\sigma(\tfrac{1}{2}) = \sigma_0$ and such that
$\hat\sigma(0)$ gives the identity of $S^5$ at the north pole of $S^6$.

The first information that we get from Lemma\,\ref{explicit-loop}
is how the actions $\alinda$ of the previous section transfer by the
boundary map $\partial$ to the cyclic subgroups of $\pi_1\Diff^+(S^{5})$
and $\pi_1\Diff^+(S^{13})$ generated by $\hat\sigma$. We evidently have
\begin{gather*}
  \widehat{\alpha \sigma \alpha^{-1}} = \alpha_0(-\hat \sigma)\alpha_0^{-1}
\end{gather*}
where $-\hat \sigma$ denotes the reverse loop, i.e.,
$(-\hat \sigma)(t)= \hat \sigma (1-t)$.
Thus, the commutation identity
$\alpha \sigma \alpha^{-1} = \sigma^{-1}$
turns into
\begin{gather*}
 \alpha_0(-\hat \sigma)\alpha_0^{-1} =(\hat \sigma)^{-1}.
\end{gather*}

\medskip

Suppose now that -- in contrast to the claim of Theorem\,\ref{5-y-13} --
the involution
$\rho_0 = \alpha_0 \sigma_0$
is conjugate to the antipodal map $\alpha_0$, i.e.,
there is a diffeomorphism $h$ of $S^5$ or $S^{13}$ such that
$\rho_0 = \alpha_0 \sigma_0 = h \alpha_0 h^{-1}$. 
Since $\alpha_0$ commutes with hyperplane reflections,
we may assume that $h$ is orientation preserving.
Solving for $\sigma_0$, we obtain 
\begin{gather*}
\sigma_0 = \alpha_0^{-1} h \alpha_0 h^{-1} \, .
\end{gather*}

Using a path $A$ in $SO(6)\subset\Diff^+(S^5)$
(resp.\ in $SO(14)\subset\Diff^+(S^{13})$)
from the identity to the antipodal map $\alpha_0$
(such a path exists for all odd dimensional spheres),
we get a path $\hat \psi$ from the identity to $\sigma_0$ given by
\begin{gather*}
\hat \psi(t) = A(t)^{-1} h A(t) h^{-1} \, .
\end{gather*}

\begin{figure}[h]
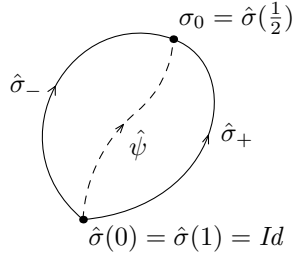
 \label{loops}
\begin{center}
 \IncludeLoopPic
 \caption{Configuration of paths in $\Diff^+(S^5)$ and $\Diff^+(S^{13})$}
\end{center}
\end{figure}

The idea is to use this path to cut $\hat \sigma$ in half and to show
that the element $[\hat \sigma] \in \pi_1 \Diff(S^{n-1})$ is the
product of a square and a correction factor contained in $N_{n-1}$,
a factorization that contradicts Lemma\,\ref{loop-two}.
In order to make this idea concrete, we decompose $\sigma$ into two paths
$\hat \sigma_+$ and $\hat \sigma_-$ defined on $[0,1]$ by
\begin{align*}
\hat \sigma_+(t) &= \hat\sigma(\tfrac{1}{2}t) \\
  \hat \sigma_-(t) &= \hat\sigma(1-\tfrac{1}{2}t)
\end{align*}
Thus, $\sigma =\hat \sigma_+ \sqcup (-\hat \sigma_-)$
where $\sqcup$ denotes juxtaposition of paths.
We concatenate each of these two paths with $-\psi$ and get loops
\begin{gather*}
 \hat \phi_\pm = \hat \sigma_\pm \sqcup (- \hat \psi) \, .
\end{gather*}
(see figure ~\ref{loops}). Clearly, we have 
\begin{gather*}
[\hat \sigma] = [\hat \sigma_+ \sqcup (-\hat \sigma_-)]
 = [\hat \phi_+ \sqcup (-\hat \phi_-)] 
 =  [\hat \phi_+ ][(-\hat \phi_-)] =  [\hat \phi_+ ][\hat \phi_-]^{-1} \, .
\end{gather*}
We claim that $[\hat \phi_+ ]$ equals $[\hat \phi_-]^{-1}$
up to some correction factor in $N_{n-1}$. This implies that
$[\hat \sigma] \equiv [\hat \phi_+]^2 \mod N_{n-1}$,
which is the contradiction that we are striving~for. 
 
Decomposing the identity
$(\hat \sigma)^{-1} = \alpha_0(-\hat \sigma)\alpha_0^{-1}$
with respect to the components of
$\hat \sigma = \hat \sigma_+ \sqcup (-\hat \sigma_-)$,
we find that on the level of paths the following holds 
\begin{gather*}
 (\hat \sigma_-)^{-1} = \alpha_0 \hat \sigma_+ \alpha_0^{-1}  \, . 
\end{gather*}
Combining these identities we obtain
\begin{gather*}
\begin{split}
  (\hat\phi_-)^{-1} &= (\hat\sigma_-)^{-1} \sqcup (-\hat\psi)^{-1}\\
  &= \alpha_0\hat\sigma_+\alpha_0^{-1} \sqcup (-\hat\psi)^{-1}\\
  &\simeq \alpha_0\bigl(\hat\sigma_+ \sqcup (-\hat\psi) \sqcup
    \hat\psi\bigr)\alpha_0^{-1} \sqcup (-\hat\psi)^{-1}\\
  &= \alpha_0\hat\phi_+\alpha_0^{-1} \sqcup \alpha_0 \hat\psi\alpha_0^{-1} 
  \sqcup (-\hat\psi)^{-1} \qquad\rel \{0,1\}.
\end{split}
\end{gather*}
Clearly, the map $(s,t) \mapsto A(s)\hat\phi_+(t)A(s)^{-1}$
provides a homotopy
$\hat\phi_+ \simeq \alpha_0\hat\phi_+\alpha_0^{-1}$.
Since $\hat\phi_+(t) = Id$ in a neighborhood of $\{0,1\}$,
the preceding map is actually a homotopy $\rel \{0,1\}$,
and so we conclude that in $\pi_1\Diff^+(S^{n-1})$
the following identity holds
\begin{gather*}
  [\hat\phi_-]^{-1} = [\hat\phi_+] \cdot [\alpha_0 \hat\psi \alpha_0^{-1}
    \sqcup (-\hat\psi)^{-1}].
\end{gather*}
Thus the proof of Theorem\,\ref{5-y-13} is completed by the following lemma:

\begin{lemma}
The path $\alpha_0\hat\psi\alpha_0^{-1} \sqcup(-\hat\psi)^{-1}$
represents an element in the normal subgroup $N_{n-1}$ of
$\pi_1 \Diff^+(S^{n-1})$.
\end{lemma}

\begin{proof}
For the purpose of this argument, we find it convenient to assume
that the path $A$ in $SO(n)\subset \Diff^+(S^{n-1})$
connecting the identity to $\alpha_0$ is the $1$-parameter subgroup
obtained by exponentiating some almost complex structure
$J \in \mathfrak{so}(n)$. On the path level, this choice yields the identity
\begin{gather*}
  (-A)^{-1} = A \alpha_0^{-1} =  \alpha_0^{-1} A\, ,
\end{gather*}
where we remind the reader that the minus sign in this equation represents
the reverse path, $-\gamma(t) = \gamma(1-t)$. 

\begin{figure}
\begin{center}
  \setlength{\unitlength}{0.01mm}
  \begin{picture}(11000,3100)
    \put(1000,0500){\line(1,0){3000}}
    \put(4000,0500){\line(0,1){2000}}
    \put(1000,0500){\line(0,1){2000}}
    \put(1000,2500){\line(1,0){3000}}
    \qbezier(2500,0500)(2450,0550)(2400,0600)
    \qbezier(2500,0500)(2450,0450)(2400,0400)
    \qbezier(1000,1500)(1050,1450)(1100,1400)
    \qbezier(1000,1500)(0950,1450)(0900,1400)
    \qbezier(4000,1500)(4050,1450)(4100,1400)
    \qbezier(4000,1500)(3950,1450)(3900,1400)
    \put(2500,0300){\makebox(0,0)[tc]{$\hat H(\,\cdot\,,0)$}}
    \put(2500,2700){\makebox(0,0)[bc]{$h\alpha_0^{-1} h^{-1} \alpha_0$}}
    \put(4300,1500){\makebox(0,0)[cl]{$\hat\psi^{-1}$}}
    \put(0800,1500){\makebox(0,0)[cr]{$\alpha_0\hat\psi\alpha_0^{-1}$}}
    \put(7000,0500){\line(1,0){3000}}
    \put(10000,0500){\line(0,1){2000}}
    \put(7000,0500){\line(0,1){2000}}
    \put(7000,2500){\line(1,0){3000}}
    \qbezier(8500,0500)(8450,0550)(8400,0600)
    \qbezier(8500,0500)(8450,0450)(8400,0400)
    \qbezier(8500,2500)(8450,2550)(8400,2600)
    \qbezier(8500,2500)(8450,2450)(8400,2400)
    \put(8500,0300){\makebox(0,0)[tc]{$\hat H(\,\cdot\,,0)
      \simeq hA^2h^{-1}A^{-2}$}}
    \put(8500,2700){\makebox(0,0)[bc]{$\alpha_0\hat\psi\alpha_0^{-1}
      \:\sqcup\: (-\hat\psi)^{-1}$\!\!}}
    \put(10300,1500){\makebox(0,0)[cl]{$Id$}}
    \put(6800,1500){\makebox(0,0)[cr]{$Id$}}
  \end{picture}
\end{center}
\caption{The homotopy $\hat H$ induces a homotopy $\rel \{0,1\}$ between
  $\alpha_0\hat\psi\alpha_0^{-1} \sqcup (-\hat\psi)^{-1}$ and
  $h A^2 h^{-1} A^{-2}$.}
\label{f:fund}
\end{figure}
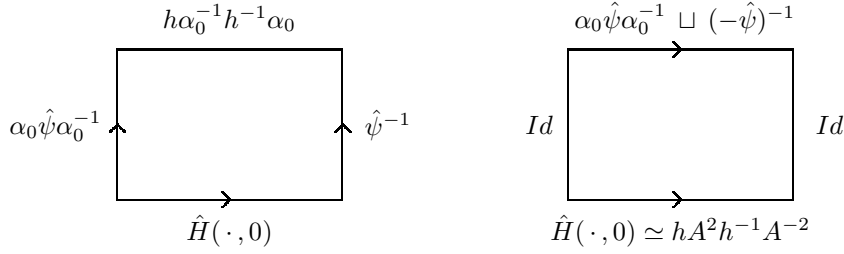

Now consider the map $\hat H: [0,1]\times [0,1] \to \Diff^+(S^{n-1})$
given by
\begin{gather*}
  \hat H(s,t) :=  \alpha_0 A(s+t-st)^{-1} h A(s+t-st)
     \alpha_0^{-1} A(1-s+ts)^{-1} h^{-1} A(1-s+ts).
\end{gather*}
Recall that $\alpha_0^2 = Id$. Thus with the help of the preceding identity,
it is easy to verify that
\begin{align*}
  \hat H(0,t) &= \alpha_0 A(t)^{-1} h A(t) \alpha_0^{-1} A(1)^{-1} h^{-1} A(1)
    = \alpha_0\hat\psi(t) \alpha_0^{-1},\\
  \hat H(1,t) &= \alpha_0 A(1)^{-1} h A(1) \alpha_0^{-1} A(t)^{-1} h^{-1} A(t)
    = \hat\psi(t)^{-1},\\
  \hat H(s,1) &= \alpha_0 A(1)^{-1} h A(1)\alpha_0^{-1} A(1)^{-1} h^{-1} A(1)
    = h\alpha_0^{-1} h^{-1} \alpha_0.
\end{align*}
Thus $\hat H$ induces a homotopy $\rel \{0,1\}$ of the concatenation
$\alpha_0\hat\psi\alpha_0^{-1} \sqcup (-\hat\psi)^{-1}$ to the path
\begin{gather*}
\begin{split}
s \mapsto \hat H(s,0)
  &= \alpha_0 A(s)^{-1} h A(s) \alpha_0^{-1} A(1-s)^{-1} h^{-1} A(1-s)\\
  &= \alpha_0 A(s)^{-1} h A(s)^2 h^{-1} A(s)^{-1} \alpha_0^{-1}\\
  &= \alpha_0 A(s)^{-1} \bigl( h A(s)^2 h^{-1} A(s)^{-2} \bigr) (\alpha_0 A(s)^{-1})^{-1}.
\end{split}
\end{gather*}
Note that $A^2$ itself is a loop in $SO(n) \subset \Diff^+(S^{n-1})$
based at the identity and so is $h A^2 h^{-1} A^{-2}$. Thus the map
\begin{gather*}
  (s,t) \mapsto \alpha_0 A(s+t-st)^{-1} \bigl(h A(s)^2 h^{-1} A(s)^{-2}\bigr)
    \bigl(\alpha_0 A(s+t-st)^{-1}\bigr)^{-1}
\end{gather*}
provides a homotopy of the path $s\mapsto \hat H(s,0)$ to the path
$h A^2 h^{-1} A^{-2}$ $\rel \{0,1\}$. Hence
\begin{gather*}
  [\alpha_0\hat\psi\alpha_0^{-1} \sqcup (-\hat\psi)^{-1}]
    = [hA^2h^{-1}] [A^2]^{-1} \in \pi_1\Diff^+(S^{n-1}).
\end{gather*}

Being the image of an element in $\pi_1 SO(n)$ under the canonical inclusion,
it is evident that $[A^2] \in  \pi_1\Diff^+(S^{n-1})$ is an element
of order at most $2$.
Conjugation by $h$ induces an automorphism of $\pi_1\Diff^+(S^{n-1})$,
and so $[hA^2h^{-1}]$ is also an element of order at most $2$.
Hence $[A^2]$ and $[hA^2h^{-1}]$ must both map into the $2$-torsion group
$T_n \subset \pi_0\Diff^+(S^n)$ under the boundary map, and therefore
\begin{gather*}
  [\alpha_0 \hat \psi \alpha_0^{-1} \sqcup (-\hat\psi)^{-1}]
   \in \partial^{-1}(T_n) = N_{n-1}
\end{gather*}
as claimed.
\end{proof}

\bigskip

\section{The geometry of Hirsch-Milnor involutions}

The fact that $\rho = \alpha\sigma$ is an exotic free involution
seems like a huge coincidence at first, stemming for the peculiar
algebraic properties of the Blakers-Massey elements $\bl$.
However, a more careful study explains why $\rho$ has these properties.
The understanding comes from the interplay between three constructions:

\begin{enumerate}

\item
The Hirsch-Milnor construction of exotic involutions of $S^5$ and $S^6$
(\cite{hirsch-milnor}), based on the Milnor exotic sphere
$\Sigma^7_{2,-1}$ (\cite{milnor}).

\item
The Gromoll-Meyer description of the Milnor exotic sphere as a quotient
$\Sigma^7_{GM}$ of the Lie group $Sp(2)$ (\cite{gromoll-meyer}).

\item
The study of the geometry of geodesics of certain metrics on the
Gromoll-Meyer exotic sphere carried out in \cite{duran}, and
\cite{duran-mendoza-rigas}, which in particular produces partial
sections and trivializations of the bundle
$S^3 \cdots Sp(2) \to \Sigma^7_{GM}$. 

\end{enumerate}

In \cite{hirsch-milnor}, Hirsch and Milnor constructed involutions of
$S^6$ and $S^5$ that are not conjugate to the antipodal map.
These involutions are constructed as follows: 
first consider the Milnor exotic sphere $\Sigma^7_{2,-1}$ (\cite{milnor}).
This sphere is a $S^3$ bundle over $S^4$ with structure group $SO(4)$,
and the bundle description is given by taking two copies of
$\mathbb R^4 \times S^3$ and identifying 
$\mathbb R^4-{0} \times S^3$ with 
$\mathbb R^4-{0} \times S^3$ via the map
\begin{gather*}
(u,v) \to (u',v') = (u/|u|^2, \frac{1}{|u|} u^2 v u^{-1}) \, .
\end{gather*}

The antipodal map in the fibers of the fibration
$S^3 \cdots \Sigma^7_{2,-1} \to S^4$
(that is, in bundle coordinates, the well-defined map 
$(u,v) \mapsto (u,-v)$) is a free involution $\tau$ on $\Sigma^7_{2,-1}$. 
This involution supports invariant spheres
$S^5 \subset S^6 \subset \Sigma^7_{2,-1}$,
described by certain subsets of the bundle coordinates, namely, 
\begin{gather*}
S^6 = \{ (u,v) \,\, \vert \,\, \Re(uv)  = 0 \} \, ,\\
S^5 = \{ (u,v) \,\, \vert \,\, \Re(uv) = \Re(v) = 0 \} \, .
\end{gather*}

Then Hirsch and Milnor show that the involutions are exotic by
showing that, if they were not exotic, then $\Sigma^7_{2,-1}$
would have even order in the cyclic group $\Gamma_7 \cong \mathbb Z_{28}$,
which contradicts the fact that $\Sigma^7_{2,-1}$ is a generator.
A proof in the same spirit shows that $\tau$ restricted to $S^5$
is also exotic.

Next in line of this exploration is the Gromoll-Meyer expression
for the sphere $\Sigma^7_{2,-1}$ as the quotient of the group $Sp(2)$
of $2\times 2$ quaternionic matrices $A$ satisfying
$A^*A = AA^* = I_{2\times 2}$ by the $S^3$-action 
\begin{gather*}
q \bigstar  
\begin{pmatrix}
a & c \\ 
b & d
\end{pmatrix}= 
\begin{pmatrix}
q & 0 \\ 
0 & q 
\end{pmatrix}
\begin{pmatrix}
a & c \\ 
b & d
\end{pmatrix}
\begin{pmatrix}
\bar q & 0 \\ 
0 & 1 
\end{pmatrix}
\begin{pmatrix}
qa\bar q & qc \\ 
qb\bar q & qd
\end{pmatrix}
\end{gather*}

The explicit identification between the Gromoll-Meyer sphere $\Sigma^7_{GM}$
and the Milnor exotic sphere $\Sigma^7_{2,-1}$ is given in
\cite{gromoll-meyer}, e.g., if $\mathcal U \subset \Sigma^7_{GM}$
is  the set of classes 
\begin{gather*}
\mathcal U \left\{
\left[
\begin{pmatrix}
a & c \\
b & d 
\end{pmatrix} 
\right]
\,\, 
\text{such that  } 
d \neq 0 
\right\}
\end{gather*}
we have 
$f:\mathcal U  \to \mathbb R^4 \times S^3$ given by 
\begin{gather*}
\left[
\begin{pmatrix}
a & c \\
b & d 
\end{pmatrix}
\right]
\mapsto
(u,v) = |d|^{-2} ( \bar c d , \bar d a d |a|^{-1} )\, ;
\end{gather*}
see \cite{gromoll-meyer} for the other charts and their inverses;
however we warn the reader that the matrices in \cite{gromoll-meyer}
are transposes of ours (we write them that way since we use the
projection onto the first column a lot).
 
From the Gromoll-Meyer formulas, we see that 

\begin{proposition} \label{from-hirsc-milnor-to-gromoll-meyer} 

The Hirsch-Milnor exotic involutions are induced by the antipodal free involution $m$ on $Sp(2)$ given by 
\begin{gather*}
\begin{pmatrix}
a & c \\ 
b & d
\end{pmatrix}
\stackrel{m}{\mapsto}
\begin{pmatrix}
-a & -c \\ 
-b & -d
\end{pmatrix}\, .
\end{gather*}
Moreover, the invariant spheres  $S^6$ (resp $S^5$) are
given by the projection of the sets $\slinda_6$ (resp. $\slinda_5$)
given by 
\begin{gather*}
\slinda_6 = 
\left\{ 
\begin{pmatrix}
a & c \\ 
b & d
\end{pmatrix} \in Sp(2) \,\, \bigg\vert \,\, \Re(a) = 0 
\right\} \, ,\\
\slinda_5 = \left\{ 
\begin{pmatrix}
a & c \\ 
b & d
\end{pmatrix} \in Sp(2) \,\, \bigg\vert \,\, \Re(a) = \Re(b) = 0 
\right\} \, .
\end{gather*}
\end{proposition}

Then we have the construction of a pointed wiedersehen metric
on $\Sigma^7_{GM}$ given in \cite{duran}. The horizontal lift
of geodesics from $\Sigma^7$ to $Sp(2)$ provides a section and
therefore a trivialization of the bundle 
$S^3 \cdots (Sp(2) \smallsetminus S^3) \to
  (\Sigma^7_{GM}\smallsetminus \{{\text{\em south pole}}\})$.

More precisely, let $N,S \in \Sigma^7_{GM}$ be given by
$N = [I_{2\times 2}], S = [-I_{2\times 2}]$.
The geodesics from $N$ are the projections of the horizontal geodesics
through the identitiy of $Sp(2)$, which are given by 
\begin{gather*}
\gamma_{(p,w)}(0) = \begin{pmatrix}
1 & 0 \\
0 & 1
\end{pmatrix},
\qquad
\gamma_{(p,w)}'(0) = \begin{pmatrix}
p & -\bar w \\
w & 0
\end{pmatrix},
\end{gather*}
where $p$ is a pure quaternion and $|p|^2 + |w|^2 = 1$. Then 
\begin{gather*}
\gamma_{(p,w)}(t) = 
\begin{pmatrix}
\cos(t) + \sin(t) p & -\sin(t)e^{tp}\bar w \\
\sin(t)w   &  \frac{w}{|w|}(\cos(t) - \sin(t)p)e^{tp} \frac{\bar w}{|w|}
\end{pmatrix}, 
\end{gather*}
in the ``generic" case $w \neq 0$. In the case $w = 0$, 
\begin{gather*}
\gamma_{(p,w)}(t) = 
\begin{pmatrix}
e^{tp} & 0\\
0  &  1
\end{pmatrix}. 
\end{gather*}

Note that the set $\slinda_6$ is the set of midpoints of the
horizontal geodesics from the identity, given by $t = \pi/2$,
and therefore the invariant $S^6 \subset \Sigma^7_{GM}$ is the
``equator" of the wiedersehen metric on $\Sigma^7_{GM}$
given by points equidistant from the north and south poles,
and the invariant sphere $S^5$ is given by $t=\pi/2$, $\Re(w) = 0$.

All we need to construct a formula for the Hirsch-Milnor
involution of $S^6$ is therefore to take advantage of the section
given by the geodesics. Let us restrict the bundle
$S^3 \cdots Sp(2) \to \Sigma^7_{GM}$ to
$S^3 \cdots \slinda_6 \to S^6$.
But this bundle is trivial, and the midpoints of the geodesics
from the north pole produce a trivialization; we have
\begin{gather*}
\psi: S^6 \times S^3 \to \slinda_6 \subset Sp(2)
\end{gather*}
given by 
\begin{gather*}
\begin{split}
\psi((p,w),q) &= 
\begin{pmatrix}
q & 0 \\ 
0 & q 
\end{pmatrix}
\gamma_{(p,w)}(\pi/2)
\begin{pmatrix}
\bar q & 0 \\ 
0 & 1 
\end{pmatrix} \\
&= \begin{pmatrix}
q & 0 \\ 
0 & q 
\end{pmatrix}
\begin{pmatrix}
  p & -e^{\frac{\pi}{2}p}\bar w \\
  w   &  -\frac{w}{|w|} pe^{\frac{\pi}{2}p} \frac{\bar w}{|w|}
\end{pmatrix}
\begin{pmatrix}
\bar q & 0 \\ 
0 & 1 
\end{pmatrix} 
\end{split}
\end{gather*}
and the inverse of $\psi$ is given by 
\begin{gather*}
\psi^{-1}
\left(
\begin{pmatrix}
a & c \\ 
b & d
\end{pmatrix} 
\right) = 
((\bar q a q, \bar q b q), q)
\intertext{where}
q
\left(
\begin{pmatrix}
a & c \\ 
b & d
\end{pmatrix} 
\right)
= -
\frac{b}{|b|} e^{-\frac{\pi}{2} a} \frac{c}{|c|} \, .
\end{gather*} 

The following is then clear:

\begin{proposition} \label{GM-S6-trivial} 
The map $\psi$ is a bundle trivialization, that is, it is an
$S^3$-equivariant map from $S^6 \times S^3 \to \slinda_6$,
where the action on the left hand side is the left multiplication
action on the $S^3$-factor and the action on the right hand side
is the Gromoll-Meyer action restricted to $\slinda_6$. 
 
\end{proposition}

The involution $m$ on $Sp(2)$ pulls back as
\begin{gather*}
\psi^{-1}\circ m \circ \psi ((p,w),q) = (\alpha\sigma^{-1} (p,w), q \,\bl(p,w))\,.
\end{gather*}
Projecting onto the first component of $S^6 \times S^3$,
we get our exotic involution $a\sigma^{-1} = \rho_{-1}$,
which is conjugate to $\rho$ by $\sigma^{-1}$.

Note that the exotic projective space $\mathbb R P_\rho^6 = S^6/\rho$
is also the quotient of $S^6 \times S^3$ under the
$\mathbb Z_2 \times S^3$-action $\star$ given by
\begin{align*}
(0, \theta) \star \bigl((p,w),q\bigr) &= \bigl((p,w), \theta q\bigr) \, ,\\
(1, \theta) \star \bigl((p,w),q\bigr) &= \bigl(\rho_{-1}((p,w)), \theta q \,\bl(p,w) \bigr) \,.
\end{align*}
since in the quotient, the $S^3$-action by $\theta$ just kills the
$S^3$-factor of $S^6 \times S^3$ just leaving
$S^6/\rho_{-1} \cong S^6/\rho$. 

\bigskip

\section{Non-cancellation phenomena in group actions} \label{noncancellation}

The usual form of non-cancellation phenomena, e.g.
\cite{hilton-mislin-roitberg-1,hilton-mislin-roitberg-2,hiltonbrazil}
is expressed by manifolds $M_1, M_2, N$ where $M_1$ and $M_2$
are not homotopy equivalent but  $M_1 \times N$ is diffeomorphic to
$M_2 \times N$ (see also \cite{barros-rigas}).
 
Here we are interested in explicit formulas for subtler
{\em differentiable} non-cancella\-tion phenomena,
exemplified for instance by the fact that for any
$7$-dimensional exotic sphere $\Sigma^7$,
$\Sigma^7 \times S^3$ is diffeomorphic to 
$S^7 \times S^3$ (see \cite{rigas, wall}).  
Here we use the exotic involution $\rho:S^6 \to S^6$
constructed above to give some differentiable non-cancellation
phenomena of {\em group actions}.

\begin{theorem}
There exist explicit actions $r_1$, $r_2$ of $\mathbb Z_2 \times S^3$ on
$X = S^6 \times S^3$ such that neither factor can be cancelled. More precisely,
\begin{itemize}
\item
the restrictions of the actions $r_1$ and $r_2$ to the subgroup
$\{0\} \times S^3$ are differentiably conjugate,
\item
the restrictions of the actions $r_1$ and $r_2$ to the subgroup
$\Z_2\times \{1_{S^3}\}$ are differentiably conjugate, 
\item
the full actions $r_1$, $r_2$ of $\mathbb Z_2 \times S^3$,
however, are {\em not} differentiably conjugate.
\end{itemize}
\end{theorem}

The construction is based on the following consideration: in addition to the
Gromoll-Meyer action, the group $S^3$ also acts freely in $Sp(2)$ as follows:
\begin{gather*}
q\bullet 
\begin{pmatrix}
a & c \\
b & d
\end{pmatrix}
= \begin{pmatrix}
a & c \\
b & d
\end{pmatrix}
\begin{pmatrix}
1 & 0 \\
0 & \bar q
\end{pmatrix},
\end{gather*}
producing a principal fibration $S^3 \cdots Sp(2) \to S^7$,
where $S^7$ is the standard $7$-sphere. In fact, the projection
of $Sp(2)$ onto $S^7$ is just $A \to$ $1^{st}$ column of $A$.

The canonical wiedersehen metric on $S^7$ produces also a partial section,
and therefore a trivialization of the bundle
$S^3 \cdots Sp(2) \smallsetminus S^3
  \longrightarrow S^7 \smallsetminus \{\text{south pole}\}$.
Let us remark that the fiber over the south pole of $S^7$ is the same
as the fiber over the south  pole of $\Sigma^7_{GM}$.

Let us list the following trivial fact as a proposition:

\begin{proposition}

The antipodal involution $m$ on $Sp(2)$ above descends under the
$\bullet$-action to the canonical involution $\alpha$ on $S^7$.  

\end{proposition}

Thus we have two different $S^3$-principal fibrations with $Sp(2)$ as total space:
\begin{gather*}
\begin{CD}
{} @. {} S^3 @. {} \\
{} @.  @V\bigstar VV @.   {} \\
S^3 @>\bullet>> Sp(2) @>>> S^7 \\
{}  @. @VVV  @.  {} \\
{} @. \Sigma^7 {}
\end{CD}
\end{gather*}

\bigskip

Note that {\em the same} involution $m$ on $Sp(2)$,
when restricted to $\slinda_6$, descends to the non-conjugate
involutions $\alpha$ and $\rho$ on $S^6$, depending on whether 
one uses the $\bullet$- or the $\bigstar$- projections, respectively.
This observation is the basis of our particular non-cancellation phenomenon. 

Let us proceed with the computation. Consider the trivialization
analogous to the one considered in the previous section, but now for
the bundle  $S^3 \cdots Sp(2) \to S^7$ restricted to
$S^3 \cdots \slinda_6 \to S^6$. In other words, consider the map 
\begin{gather*}
  \phi: S^6 \times S^3 \to \slinda_6 \subset Sp(2),\\
\phi((p,w),q) \,=\, 
\gamma_{(p,w)}(\pi/2)
\begin{pmatrix}
1 & 0 \\ 
0 & \bar q 
\end{pmatrix}
\,=\,
\begin{pmatrix}
  p & -e^{\frac{\pi}{2}p}\bar w \\
  w   &  -\frac{w}{|w|} pe^{\frac{\pi}{2}p} \frac{\bar w}{|w|}
\end{pmatrix}
\begin{pmatrix}
1 & 0 \\ 
0 & \bar q 
\end{pmatrix}.
\end{gather*}
The inverse of $\phi$ is given by 
\begin{gather*}
\phi^{-1}
\left(
\begin{pmatrix}
a & c \\ 
b & d
\end{pmatrix}
\right) = 
((a,b), q), \quad \text{where}\quad
q
\left(
\begin{pmatrix}
a & c \\ 
b & d
\end{pmatrix} 
\right) = 
-
\frac{{\bar c}}{|c|} e^{\frac{\pi}{2} a} \frac{\bar b}{|b|} \, .
\end{gather*}

Similarly, we have 

\begin{proposition} \label{standard-is-trivial}
The map $\phi$ is a bundle trivialization, that is,
it is an $S^3$-equivariant map from $S^6 \times S^3 \to \slinda_6$,
where the action on the left hand side is the left multiplication
action on the $S^3$-factor and the action on the right hand side
is the standard $\bullet$-action above action restricted to $\slinda_6$. 
\end{proposition}

Then a computation shows that 
\begin{gather*}
\phi^{-1}\circ m \circ \phi ((p,w),q)
  = (\alpha(p,w), q \,\overline{\bl(p,w)})\, ,
\end{gather*}
and projecting onto the first component of $S^6 \times S^3$
we get the standard antipodal involution on $S^6$. 

The standard projective space $\mathbb R P^6 = S^6 / \alpha$
is then the quotient of the
$\mathbb Z_2 \times S^3$-action $\circdot$ on $S^6 \times S^3$
given by 
\begin{align*}
(0, \theta) \circdot \bigl((p,w),q\bigr) &= \bigl((p,w), \theta q\bigr) \, ,\\
(1, \theta) \circdot \bigl((p,w),q\bigr) &= 
  \bigl(\alpha(p,w), \theta q \,\overline{\bl(p,w)} \bigr) \,.
\end{align*}

Let $r_1, r_2$ now be, respectively, the $\star$- and $\circdot$-actions
defined above. Then clearly $r_1$ and $r_2$ are not conjugate,
since a conjugacy between $r_1$ and $r_2$ would imply that
$\mathbb R P^6_\rho$ and $\mathbb R P^6_\alpha$
are diffeomorphic, and we know they are not. 

Observe that the restrictions of $r_1$ and $r_2$ to $\{0\} \times S^3$
coincide and so they are trivially differentiably conjugate.
On the other hand, the restrictions of $r_1$ and $r_2$ to
$\Z_2 \times \{1_{S^3}\}$ are given by
\begin{align*}
  (1,1_{S^3}) \star \bigl((p,w),q\bigr)
  &= \bigl(\rho_{-1}(p,w), q\, b(p,w)\bigr),\\
  (1,1_{S^3}) \circdot \bigl((p,w),q\bigr)
  &= \bigl(\alpha(p,w), q\, \overline{b(p,w)}\bigr).
\end{align*}

These two actions are equivalent; they are really the
same involution $m$ on $Sp(2)$ restricted to the set
$\mathcal S_6 \cong S^6 \subset Sp(2)$.
The fact that they look different comes from the different
trivializations of the bundle
$S^3 \cdots \mathcal S_6 \to S^6$.
In fact, the diffeomorphism conjugating them is given by 
\begin{gather*}
F((p,w),q) = ((qp\bar q, q w \bar q , \bar q) \, .
\end{gather*}

Note that the conjugating diffeomorphism $F$ is also an involution!

\end{document}